\documentclass[11pt,twoside]{article}
\usepackage{amsfonts,amssymb}

\def\Bel{\hbox{\vrule\vbox to 7pt
{\hrule width 6pt
\vfill\hrule}\vrule}} 
\def\Bbb{\mathbb}
\def\qed{\space\Bel}
\def\goth{\mathfrak}
\title{Trace class groups: the case of semi-direct products}
\newtheorem{stelling}{Theorem}
\newtheorem{hulpst}{Lemma}
\newtheorem{cor}{Corollary}

\newtheorem{proposition}{Proposition}
\author{Gerrit van Dijk}
\date{}

\begin{document}
\maketitle
\par\noindent
{\bf Abstract}.
A Lie group $G$ is called a trace class group if for every irreducible unitary representation $\pi$ of $G$ and every $C^\infty$ function $f$ with compact support the operator $\pi (f)$ is of trace class. In this paper we extend the study of trace class groups, begun in a previous paper, to special families of semi-direct products. For the case of a semisimple Lie group $G$ acting on its Lie algebra $\goth g$ by means of the adjoint representation we obtain a nice criteriom in order that $\goth g\ltimes G$ is a trace class group.
\vskip1cm
\noindent
\centerline{{\bf Introduction}}
\vskip.3cm\noindent
In this paper we continue the study of trace class groups begun in [DD15]. An irreducible unitary representation $\pi$ is said to be of trace class if for every $C^\infty$ function $f$ with compact support the operator $\pi (f)$ is of trace class. A group is said to be of trace class, or briefly, a trace class group, if every irreducible unitary representation is of trace class. Well-known examples of such groups are reductive Lie groups and unipotent Lie groups. It is known that the direct product of two trace class groups is again a trace class group. This needs however not to be true for semi-direct products. Since in general each Lie group is a semi-direct product of a reductive and a unipotent Lie group, we will again take up the study of semi-direct products in this paper and start our investigations about a criterion for a semi-direct product in order to be a trace class group.  We shall, for the time being, restrict to groups for which the unipotent part (the unipotent radical) is abelian. Several new examples of such groups are studied. One of the highlights is the theorem that the semi-direct product of a semisimple Lie group and its Lie algebra is a trace class group if and only if the group is compact. A nice generalization of this theorem would be provided by the case of a semisimple Lie group $G$ acting on a real finite-dimensional vector space $V$ by linear transformations and considering the semi-direct product of $V$ and $G$. A proof is not yet available however at this time.

\section{Levi decomposition}
We begin with recalling the Chevalley or Levi decompostion of an algebraic linear group over $\Bbb R$ or $\Bbb C$. 
\par
Let $G$ be such a group. The {\it unipotent radical} ${\rm Rad}_u (G)$ of $G$ is the largest {\it normal} subgroup of $G$ consisting of unipotent elements. The unipotent  radical is a connected algebraic subgroup. If $G$ is reductive then its unipotent radical consists of the identity element and conversily. The unipotent radical is itself a unipotent group.  
\setcounter{stelling}{1}
\begin{proposition} There exists a reductive algebraic subgroup $H$ of $G$ such that $G={\rm Rad}_u (G)\ltimes H$. 
\end{proposition}
This is the Chevalley or Levi decomposition of $G$. For a proof see [OV90], Ch.1, 6.5. We conclude again that in order to find a criterion for a group to be of trace class, the study of semi-direct products is crucial.
\section{Irreducible unitary representations of semi-direct products}
We start with the more general case of semi-direct products of Lie groups. 
\par\noindent
Let $G$ be a Lie group being the semi-direct product $G=N\ltimes H$ of a closed normal  subgroup $N$ and a closed subgroup $H$ and assume that $N$ is {\it abelian}. The group $H$ acts continuously on $N$ by $h\mapsto \alpha_h$ where $\alpha_h(n)=hnh^{-1}\ (n\in N)$. Denote by $\widehat N$ the dual group of $N$ and define an action of $H$ on $\widehat N$ by $\alpha^\ast_h(\chi)=\chi \circ \alpha_{h^{-1}}$.
\par\noindent
The {\it orbit} $O_\chi$ of $\chi\in\widehat N$ is the set of all elements $\alpha_h^\ast(\chi)$, for $h$ running through $H$. Define the {\it stabilizer} $H_\chi$ of $\chi$ as the set $H_\chi =\{h\in H\, |\, \alpha^\ast_h\chi=\chi\}$. Then $O_\chi\simeq H/H_\chi$.
\par\noindent
 We shall describe a wide class of irreducible unitary representations of $G$. One has the following proposition.
\setcounter{stelling}{1}
\begin{proposition} Let $G=N\ltimes H$ be a semi-direct product of $N$ and $H$ and assume that $N$ is abelian, $H$ unimodular and the action of $H$ on $N$ is unimodular. Let $\chi\in\widehat N$ and assume that the orbit $O_\chi$ carries a $H$-invariant measure. Then one has:
\par 
(i) If $\rho$ is an irreducible unitary representation of $H_\chi$, then $\chi\otimes \rho$ is an irreducible unitary representation of $N\ltimes H_\chi$. 
\par
(ii) The representation induced by $\chi\otimes\rho$ is an irreducible unitary representation of $G$.
\par
(iii) If $O_{\chi_1}=O_{\chi_2}$, and $\rho_1$ is an irreducible unitary representation of $H_{\chi_1}$, then there exists an irreducible unitary representation $\rho_2$ of $H_{\chi_2}$ such that the representations induced by $\chi_1\otimes\rho_1$ and $\chi_2\otimes\rho_2$ respectively are equivalent.
\par\noindent 
If, in addition, any orbit $O_\chi$  carries an $H$-invariant measure, then one has:
\par
(iv) If $C$ is a set in $\widehat N$ meeting each orbit only once, then the representations induced by $\chi\otimes\rho\ (\chi\in C)$ are mutually inequivalent.
\par
(v) If one can choose $C$ as a Borel subset, each irreducible unitary representation of $G$ is equivalent with an induced representation as obtained above.
\end{proposition}
\noindent
For a proof we refer to [Mack52] and [Mack68].
\par\noindent
If ($v$) holds, then one speaks of a {\it Mackey-regular semi-direct product}.  
\par\noindent
We remark that the requirement that any orbit $O_\chi$ admits a $H$-invariant measure, restricts the class of semi-direct products we consider. Indeed there exist semi-direct products where some orbits do not carry such a measure. We refer to [vD83], Section 4. 
\vskip.1cm\noindent
The next step is to determine whether the above induced representations are of trace class and, if so, to  establish a suitable formula for the trace. 
\section{A formula for the trace}
Let again $G=N\ltimes H$ with $N$ abelian, $H$ unimodular and the action of $H$ on $N$ unimodular. Then $G$ is unimodular too.
\par\noindent
Fix $\chi_0\in\widehat N$ and denote by $H_0$ the stabilizer of $\chi_0$ in $H$. Assume that $H/H_0$ carries a $H$-invariant measure. Let $\rho$ be an irreducible unitary representation of $H_0$ and define $\pi$ as the irreducible unitary representation of $G$  induced by $\chi_0\otimes\rho$. Then one has:
\setcounter{stelling}{0}
\begin{stelling}
For any function $\varphi\in C_c^\infty (G)$ one has
\setcounter{equation}{0}
\begin{equation}
trace\,\pi (\varphi)=trace\, \int_{H/H_0}\int_N\int_{H_0}\,\varphi (n, hh_0h^{-1})\,\chi_0(h^{-1}nh)\, \rho(h_0)\, dh_0dn d{\dot h}.
\end{equation}
The left hand side is finite if and only if the right-hand side is.
\end{stelling}
The proof is for the greater part in [Kir76], \S 13. Application of Mercer's theorem and the observation that any $\varphi$ is a linear combination of function of the $\psi *\widetilde\psi$ completes the proof.
\section{Application}
We give an alternative proof of [DD15], Proposition 1.13: 
\setcounter{stelling}{1}
\begin{proposition} The group ${\Bbb R}^n\ltimes {\rm SL}(n,\Bbb R )$ is not a trace class group.
\end{proposition}
{\bf Proof}. Set $N=\Bbb R^n$ and $H={\rm SL}\, (n,\Bbb R )$. Let $\chi_0 (x)={\rm e}^{-2\pi i (x\, |\, e_1)}\ (x\in {\Bbb R}^n)$, where $e_1=(1,0,\dots,0)\in{\Bbb R}^n$ and let $\rho=id$, the identity representation of $H_0$. We shall show that $\pi={\rm Ind}_{N\ltimes H_0}^G\, \chi_0\otimes\rho$ has no finite trace. Observe that $H/H_0\simeq {\Bbb R}^n\backslash\{0\}$  carries a $H$-invariant measure, namely the ususal Lebesgue measure. We have to consider the right-hand side of equation (1). We get
\begin{equation}
\int_{H/H_0}\int_H\int_{H_0}\,\varphi (n, hh_0h^{-1})\, \chi_0(h^{-1}nh)\, dh_0dn d{\dot h}.
\end{equation}
Clearly $H_0$ consists of the matrices of the form 
$$\left (\begin{array}{cc}1&0\\x&B\end{array}\right)$$
 where $x\in {\Bbb R}^{n-1}$ and $B\in {\rm SL}\, (n-1,\Bbb R^n)$. Setting 
$$M_0:\{m=\left (\begin{array}{cc}1&0\\0&B\end{array}\right)\}\quad {\rm and}\quad U:\{u=u(x)=\left (\begin{array}{cc}1&0\\x&I\end{array}\right)\}$$
we may write $H_0=M_0U$. Both $M_0$ and $U$ are unimodular. Select Haar measues $dm$ on $M_0$ and $dx$ on $U$. Then $dh_0=dmdx$. Let $A$ be the subgroup of ${\rm SL}\, (n,\Bbb R)$ consisting of the matrices of the form
$$\left (\begin{array}{cc}\lambda^{-1}&0\\0&\mu\, I\end{array}\right)$$
with $\lambda\not=0,\,\mu\not=0,\, \mu^{n-1}=\lambda$. Let $K={\rm SO}\, (n,\Bbb R)$. One has $H=KAM_0U$ and $dh=\vert\lambda\vert^{n-1}\, dkd\lambda dmdx$, hence $d{\dot h}=\vert\lambda\vert^{n-1}dk d\lambda$. Moreover, it is easily checked that $M_0$ and $A$ commute. For a function $\varphi$ of the form $\varphi(n,h)=\varphi_1(n)\,\varphi_2(h)\ (n\in N, h\in H)$ where $\varphi_1$ and $\varphi_2$ are assumed to be $K$-invariant, we obtain from (2):
$$\int \,\widehat{\varphi_1}(\lambda e_1)\cdot$$
$$\{\int \,\varphi_2 (m\, u(\lambda\mu\,x_1,\lambda\mu\, x_2,\ldots,\lambda\mu\, x_{n-1}))\, dm\,dx_1dx_2\ldots dx_{n-1}\}\,\vert\lambda\vert^{n-1}d\lambda$$
$$=\int_{-\infty}^\infty \widehat{\varphi_1}(\lambda e_1).\{\int_{MU}\varphi_2(mu(x)) dm\,dx\}\, {d\lambda\over \vert\lambda\vert}.$$
This expression is clearly divergent for suitable $\varphi_1$ en $\varphi_2$. Therefore the ${\rm trace}$ of $\pi$ does not exist.\qed
\vskip.1cm\noindent
If $n=2$, so if $G={\Bbb R}^2\ltimes {\rm SL}\, (2,\Bbb R )$, the representation $\pi$ considered above is the only one without a trace.
\par\noindent
The method of proof of Propososition 4.1 can be applied in several other situations. See Section 5.
\section{More examples}
For details about the structure of the groups considered below, see [Fa64].
\par\noindent
Let $\Bbb F$ stand for the field $\Bbb R$ or $\Bbb C$ and let $H$ be one of the groups ${\rm O}(p,q)$ or ${\rm U}(p,q) $, abbreviated by $U(p,q;\Bbb F)$. Set $n=p+q$. Denote by $E=E_n$ the space ${\Bbb F}^n$. The group $H$ acts on $E$ in the standard way. Let us consider the semi-direct product $G=E\ltimes H$. If the elements of $G$ are written as pairs $g=(v,h)\ (v\in E,\, h\in H)$ then the product in $G$ is given by
$$(v,h)\, (v',h')=(v +h\cdot v',\, hh')\quad (v,v'\in E; h,h'\in H).$$
Notice that the group $G=E\ltimes H$ is unimodular, that $H$ is unimodular and acts on $E$ in a unimodular way.
\setcounter{stelling}{1}
\begin{proposition}
The group ${\Bbb F}^n\ltimes {\rm U}(p,q)$ is not a trace class group unless $p=0$ or $q=0$.
\end{proposition}
{\bf Proof}. Define the goups $N$ and $M$ as follows:
\vskip.1cm\noindent
$$M=\{\left(\begin{array}{ccc}u&0&0\\0&v&0\\0&0&u\end{array}\right ) : u\in{\Bbb F}, \vert u\vert =1,\, v\in U(p-1,q-1,\Bbb F)\},$$
$$N=\{ n(z,w)=\left ( \begin{array}{ccc}1+w- {1\over 2}[z,z]&z^\ast&-w+{1\over 2}[z,z]\\z&I&-z\\w-{1\over 2}[z,z]&z^\ast&1-w+{1\over 2}[z,z]\end{array}\right )\}$$
where $w\in\Bbb F,\, w+\overline w=0,\, z^\ast =-\ ^t{\overline z}I_{p-1,q-1}$ and 
if 
\par
$z=\left (\begin{array}{c}z_1\\\vdots\\z_{n-2}\end{array}\right ), \, z'=\left (\begin{array}{c}z'_1\\\vdots\\z'_{n-2}\end{array}\right )$, then
\vskip.1cm
$[z,z']=\overline z'_1 z_1 +\cdots +\overline z_{p-1}' z_{p-1} - \overline z_{p}' z_{p} -\cdots - \overline z'_{n-2} z_{n-2}\quad (n=p+q).$
\vskip.2cm
One has
$$n(z,w) \cdot n(z',w') = n(z + z',w + w' + Im[z, z']).$$
The group $M$ acts on $N$:
$$m\cdot n(z,w)\cdot m^{-1}= n(u^{-1}vz, w)$$
where $m=\left(\begin{array}{ccc}u&0&0\\0&v&0\\0&0&u\end{array}\right )\in M$. Then $MN$ is clearly a (unimodular) Lie group, being the semi-direct product of $M$ and $N$.
 Let $\xi^0=(1,0,\ldots ,0,1)\in E$. The stabilizer of $\xi^0$ in $H$ is the group $H_0=M_0N$ where $M_0$ is the subgroup of $M$ with $u=1$, so $M_0\simeq {\rm U}(p-1,q-1)$. Let $A$ be the group of matrices of the form
$$a_t=\left(\begin{array}{ccc}{\rm cosh}\, t&0&{\rm sinh}\, t\\0&I_{p-1,q-1}&0\\{\rm sinh}\,t&0&{\rm cosh}\, t\end{array}\right )$$
where $t\in\Bbb R$. Then $P=MAN$ is a parabolic subgroup of $H$, being the stabilizer of the line ${\Bbb F}\,\xi_0$. The group $A$ normalizes $N$:
$$a_t\, n(z,w)\, a_{-t}=n({\rm e}^{t}z ,{\rm e}^{2t}w),$$
and commutes with $M$.
Let $K={\rm U}(p)\times {\rm U}(q)$. Then one has $H=KAMN=KAM_0N$. Select Haar measures $dk$ on $K $, $dm$ on $M_0$, $dt$ on $A$ and $dn=dzdw$ on $N$. Then $dh={\rm e}^{2\rho t}\,dkdtdmdn$, where $2\rho=n-2\ ({\Bbb F}={\Bbb R}),\,2\rho=2n-2\ ({\Bbb F}={\Bbb C})$.
\par\noindent
Define the character $\chi_0$ of $E$ as follows. Set
$$[v,v']=\overline v'_1 v_1 +\cdots +\overline v_{p}' v_{p} - \overline v_{p+1}' v_{p+1} -\cdots - \overline v'_{n} v_{n}$$
for $v,v'\in E$. Then $\chi_0(v)={\rm e}^{-2\pi i{\rm Re}\,[v,\,\xi^0]}\ (v\in E)$. Consider again $\pi={\rm Ind}_{E\ltimes H_0}^G\chi_0\otimes id$. Similarly to the proof of Proposition 4.1 we again meet a divergent integral for ${\rm trace}\,\pi$.\qed
\vskip.2cm\noindent
Other examples are ${\Bbb C}^n\ltimes O(n,{\Bbb C})$ and ${\Bbb C}^n\ltimes SL(n,{\Bbb C})$.\section{Semi-direct product of abelian groups}
The following result is in the thesis of Klamer [Kla79].
\par\noindent
Let $G=N\ltimes H$ be a semi-direct product of the locally compact abelian groups $N$ and $H$, and assume that the Haar measure on $N$ is $H$-invariant. Then $G$ is a unimodular group. 
\par
Let $\chi_0\in\widehat{N}$ and set $H_0$ for the stabilizer of $\chi_o$ in $H$. Select $\rho\in \widehat{H}_0$ and define the representation $\pi=\pi_{\chi_0,\rho}$ as 
$$\pi={\rm Ind}_{N\ltimes H_0}^ G\, \chi_0\otimes \rho.$$
We know that $\pi$ is irreducible and unitary. Denote by $O_{\chi_0}$ the $H$-orbit of $\chi_0$ in $\widehat{N}$. Since $O_{\chi_0}\simeq H/H_0$, the orbit obviously carries a $H$-invariant measure. 
\setcounter{stelling}{0}
\begin{stelling}
The representation $\pi_{\chi_0,\rho}$ has a trace if and only if the $H$-invariant measure on $O_{\chi_0}$ naturally extends to a tempered $H$-invariant measure on $\widehat{N}$.
\end{stelling}
{\bf Proof}. Denote by $\mu$ the $H$-invariant meassure on $O_{\chi_0}$. If $\mu$ extends to a tempered $H$-invariant measure on $\widehat N$, then clearly ${\rm trace}\, \pi$ exists. Indeed, apply formula (1) and observe that
$$n\mapsto\int_{H_0} \varphi (n,h_0)\,\rho(h_0)\, dh_0$$
is in $C_c^\infty (N)$, and therefore its Fourier transform is in ${\cal S}(\widehat N)$.
\par
For the converse direction assume that ${\rm trace}\,\pi (\varphi)$ exists for all $\varphi\in C_c^\infty (G)$. Taking $\varphi$ of the form $\varphi =\varphi_1\otimes\varphi_2$ with $\varphi_1\in C_c^\infty (N),\,\varphi_2\in C_c^\infty (H)$, we see from the formula for the trace that
$\int_{O_{\chi_0}}\widehat{\varphi}_1(\chi)\, d\mu(\chi)$ exists for all $\varphi_1\in C_c^\infty (N)$, and, in addition, that this expression defines a distribution $T$ on $N$, since ${\rm trace}\,\pi$ is a distribution.
\par
Observe that for all $\chi\in\widehat N$ there exists a function $\varphi_1\in C_c^\infty (N)$ with $\widehat{\varphi}_1(\chi)\not=0$. It follows that $\mu$ is a locally finite measure on $\widehat N$, hence a (Radon) measure on $\widehat N$. The distribution $T$ is clearly positive-definite, hence tempered [Sch78,Th\'eor\`eme VII, p. 242].  Thus $\mu$ is a tempered measure on $\widehat N$, since $\varphi_1\mapsto\widehat{\varphi}_1$ is an isomorphism of ${\cal S}(\widehat{N})$ onto ${\cal S}(N)$. \qed
\vskip.2cm\noindent 
This interesting theorem enables us to find examples of non-trace class groups in an easy way. Take for instance $G={\Bbb R}^2\ltimes O(1,1)$. 
\par
On the positive side we have the following examples of semi-direct products of abelian groups.
\vskip.2cm\noindent
1. {\bf Heisenberg groups}
\vskip.2cm\noindent
Let
$$G=\left\{\left(\begin{array}{ccc}1&a_1\ldots a_{n-1}&b_1\\0&1&\vdots \\0&0&b_n\\0&0&1\end{array}\right):\, a_i,b_j\in\Bbb R\right\}$$
and denote by $N$ and $H$ the subgroups defined by $a_1=\ldots = a_{n-1}=0$ and $b_1=\ldots =b_n=0$ respectively. Then both $N$ and $H$ are closed abelian subgroups of $G$. $N$ is a normal subgroup isomophic to ${\Bbb R}^n$ and $H$ is isomorphic to ${\Bbb R}^{n-1}$. Moreover $G=NH$ and $N\cap H=\{1\}$. $G$ is a semi-direct product of $N$ and $H$. It turns out that it is even a regular semi-direct product.
\vskip.2cm\noindent
2. {\bf The group of $4\times 4$ upper triangular unipotent matrices}
\vskip.2cm\noindent
Let $G$ denote the group of matrices
 $$G=\left\{\left(\begin{array}{cccc}1&x&a&b\\0&1&c&d\\0&0&1&y\\0&0&0&1\end{array}\right):\, x,y,a,b,c,d\in\Bbb R\right\}$$
and set 
$$N=\left\{\left(\begin{array}{cccc}1&0&a&b\\0&1&c&d\\0&0&1&0\\0&0&0&1\end{array}\right):\, a,b,c,d\in\Bbb R\right\}.$$
Then $N$ is a closed abelian normal subgroup of $G$ isosmorphic to ${\Bbb R}^4$. Let
$$H=\left\{\left(\begin{array}{cccc}1&x&0&0\\0&1&0&0\\0&0&1&y\\0&0&0&1\end{array}\right):\, x,y\in\Bbb R\right\}.$$
Then $H$ is a closed abelian subgroup isomorphic to ${\Bbb R}^2$. One has $G=NH$, $N\cap H=\{1\}$. $G$ is naturally a semi-direct product of $N$ and $H$. It turns out that it is even a regular semi-direct product.
In both examples $G$ is clearly a trace-class group, hence the $H$-invariant measures on the $H$-orbits in $\widehat{N}$ are tempered. This fact can of course also be checked directly.
\section{Compact stabilizers}
Let $G=N\ltimes H$, $N$ abelian, $H$ unimodular and the action of $H$ on $N$ also unimodular. Suppose  that for  $\chi_0 \in\widehat{N}$ the stabilizer of $\chi_0$ in $H$ is {\bf compact}. Select $\rho\in\widehat{H}_0$ and define $\pi=\pi_{\chi_0,\rho}$ as above. 
\setcounter{stelling}{0}
\begin{stelling}
The representation $\pi$ has a trace if the $H$-invariant measure on $O_{\chi_0}$ naturally extends to a tempered $H$-invariant measure on $\widehat N$.
\end{stelling}
{\bf Proof}. Apply formula (1) and observe that for all $\psi\in C_c^\infty (G)$
$${\rm trace}\, \int_{H_0}\,\psi(n,hh_0h^{-1})\,\rho(h_0)\, dh_0$$
is a bounded continuous function of $n$ and ${\dot h}$, of compact support with respect to $n$. In case $N=\Bbb R^n$ take for instance $\psi=(1+\Delta_N)^l\,\varphi$ where $\Delta$ is the usual Laplace operator on $\Bbb R^n$, $\varphi\in C_c^\infty(G)$ and $l$ a sufficiently large integer.\qed
\vskip.1cm\noindent
As application consider the group ${\Bbb R}^{n+1}\ltimes O(1,n)$ and set for instance $\chi_0(x)={\rm e}^{-2\pi i x_0}$ for all $x=(x_0,x_1,\ldots,x_n)\in{\Bbb R}^{n+1}$. 
\par\noindent
Another application is provided by the semi-direct products $H\ltimes \frak h$ where $H$ is a connected semi-simple Lie group with finite center and $\frak h$ its Lie algebra. $H$ acts on $\frak h$ by the adjoint representation. It is well-known that any regular $H$-orbit on $\frak h$ admits an invariant measure, that can be extended to $\frak h$ as a tempered Radon measure. Some groups $H$ admit a compact Cartan subgroup $A$, namely if $H$ admits a discrete series, for example ${\rm SL}(2,\Bbb R)$. Let $\frak a$ be the Lie algebra of $A$ and select a regular element $X$ in $\frak a$. Then the stabilizer of $X$ is precisely $A$, so compact. 
\setcounter{stelling}{2}
\begin{cor} Let $G=N\ltimes H$, $N$ abelian and $H$ compact. Then every $\pi=\pi_{\chi_0,\rho}$ has a trace.
\end{cor}
{\bf Proof}. This is obvious by Theorem 7.1 and the observation that any measure with compact support is tempered.\qed
\section{The semi-direct product ${\Bbb R}^3\ltimes {\rm SO}_0(1,2)$}
Let us write $N={\Bbb R}^3$, $H= {\rm SO}_0(1,2)$ and $G=N\ltimes H$. The $H$-orbits on ${\Bbb R}^3$ are given as follows. Set
$$[x,\, y]=x_0y_0 - x_1y_1 - x_2y_2\quad  (x,y\in{\Bbb R}^3).$$
Then the non-trivial orbits are given by $[x,\, x]=\alpha^2,\  x_0>0,\ [x,\, x]=\alpha^2,\  x_0<0,\ [x,\, x]=0, x_0>0,\ [x.\, x]=0,\, x_0<0,\ {\rm and}\ [x,\, x]=-\alpha^2$, where throughout $\alpha$ is a strictly positive real number. Call the latter orbit $\cal O_\alpha$. The first two orbits give compact stabilizers and it is well-known that the $H$-invariant measures naturally extend to tempered invariant measures on $\widehat{N}$. The third and fourth orbit lead to representations of $G$ without trace. We will consider the case of the orbit ${\cal O}_\alpha$. The stabilzer of $\alpha\,e_1$ is equal to $H_0=A$ where  $A$ is the group of matrices of the form
$$a_t=\left(\begin{array}{ccc}{\rm cosh}\, t&0&{\rm sinh}\, t\\0&1&0\\{\rm sinh}\,t&0&{\rm cosh}\, t\end{array}\right )$$
with $t\in\Bbb R$. 
\par
The irreducible unitary representations of $H_0$ are given by 
$$\rho_s:\, a_t\mapsto {\rm e}^{-2\pi i st}\quad (t\in{\Bbb R}).$$
Define
	$$\chi_\alpha (x)={\rm e}^{-2\pi 1 \, [x,\, e_2]\alpha}\quad (x\in{\Bbb R}^3,\, \alpha>0).$$
The stabilizer of $\chi_\alpha$ is $H_0$ again. Set $\pi_{s,\alpha}={\rm Ind}_{{\Bbb R}^3\ltimes H_0}^G\, \chi_\alpha \otimes \rho_s$. We shall show that ${\rm trace}\, \pi_{s,\alpha}$ exists.
\par
We start with formula (6.1) from Section 6.3. For any $\varphi\in C_c^\infty (G))$ one has
$${\rm trace}\,\pi_{s,\alpha)}(\varphi)=\int_{H/H_0}\int_{{\Bbb R}^3}\int_{\Bbb R}\,\varphi(x;ha_th^{-1})\, \chi_\alpha (h^{-1}\cdot x)\, \rho_s(a_t)\, dtdx d{\dot h}.$$
Write $H=KUA$ (Iwasawa decomposition) with $K\simeq {\rm SO}(2)$ and 
$$U=\left\{u_z=\left(\begin{array}{ccc}1+{z^2\over 2}&z&-{z^2\over 2}\\z&1&-z\\{z^2\over 2}&z&1-{z^2\over 2}\end{array}\right ):\, z\in\Bbb R\right\}.$$
Then $dh=dkdzdt$. We may always assume that $\varphi$ is ${\rm Ad}(K)$-invariant. Then we obtain:
$${\rm trace}\,\pi_{s,\alpha}(\varphi)=\int_{\Bbb R}\int_{{\Bbb R}^3}\int_{\Bbb R}\,\varphi(x;\, u_za_tu_{-z})\,\chi_{\alpha}(u_{-z}\cdot x)\,\rho_s(a_t)\, dtdxdz.$$
Let $D_x=\displaystyle {1\over{2\pi^2\alpha^2}}\left ( {\partial^2\over\partial x_1^2}+{\partial^2\over\partial x_2^2}\right )$, an ${\rm Ad}(K)$-invariant differential operator. Then we get:
$${\rm trace}\,\pi_{s,\alpha}(\varphi)$$ $$=\int_{\Bbb R}\int_{{\Bbb R}^3}\int_{\Bbb R}\,{1\over 1+z^2}\, D_x\varphi(x;\, u_za_tu_{-z}a_{-t}\cdot a_t\,\chi_{\alpha}(u_{-z}\cdot x)\,\rho_s(a_t)\, dtdxdz.$$
Since $u_za_tu_{-z}a_{-t}\in U$, we see that $a_t$, hence $t$, varies in a compact set. This is a crucial point. Hence 
$$\vert {\rm trace}\,\pi_{s,\alpha}(\varphi)\vert\leq C\cdot{\rm max}_{s,h}\, \left\vert D_x\,\varphi(x;h)\right\vert \cdot \int_{\Bbb R}\, {dz\over 1+z^2}<\infty.$$
This shows that ${\rm trace}\, \pi_{s,\alpha}$ exists.
\section{A special family of semi-direct products}
\setcounter{stelling}{0}
In this section we shall prove the following theorem:
\begin{stelling}
Let $G$ be a connected real semisimple Lie group with finite center and denote by $\goth g$ its Lie algebra. Let the group $G$ act on $\goth g$ by means of the adjoint representation. Then the semi-direct product ${\goth g}\ltimes G$ is a trace class group if and only if $G$ is compact.
\end{stelling}
{\bf Proof}. It is well-known that $G$ is unimodular and acts on $\goth g$ in a unimodular way. Furthermore, every $G$-orbit on the dual of $\goth g$ carries a $G$-invariant measure. Let us assume that $G$ is non-compact. Choose a Cartan involution $\theta$ of $\goth g$ and let ${\goth g}={\goth k} + {\goth p}$ be its decomposition into $\pm 1$-eigenspaces. Denote by $K$ the compact subgroup with Lie algebra $\goth k$. Select a maximal abelian subspace $\goth a$ of $\goth p$ and let $\Sigma$ denote the set of roots of $({\goth g},\,{\goth a})$. Then $\Sigma$ is a root system (with multiplicities). Let $\Delta$ be a set of simple roots and $\Sigma^+$ the set of positive roots with respect to $\Delta$. Set $\goth u$ for the space spanned by the positive root vectors. Then ${\goth g}={\goth k}+{\goth a}+{\goth u}$ and similarly $G=KAU$, the Iwasawa decomposition of $G$.
Let $\beta$ be a maximal positive root, i.e. $\alpha +\beta$ is not a root for any $\alpha\in{\Sigma^+}$. Denote by $X_0\not= 0$ a root vector for $\beta$ and by $G_0$ the stabilizer of $X_0$ in $G$. Let $\langle\ ,\ \rangle$ denote the Killing form of $\goth g$. Define $H_0\in{\goth a}$ by
$$\beta (H)=\langle H,H_0\rangle \quad (H\in{\goth a})$$
and set $A_1={\rm exp}\, {\Bbb R}H_0$. Define the subgroup $P$ as the stabilizer of the half-line ${\Bbb  R}^\ast_+\, X_0$. Then $P$ is equal to the semi-direct product $P=A_1G_0$, $A_1\cap G_0=\{1\},\, a_1G_oa_1^{-1}=G_0$ for all $a_1\in A_1$. Let $da_1$ and $dg_0$ denote left Haar measures on $A_1$ and $G_0$ respectively. Then $dp=da_1dg_0$ is a left-invariant measure on $P$. Since every $G$-orbit on $\goth g$ carries a $G$-invariant measure, $dg_0$ is also right invariant. So $d_rp=\Delta(a_1)\, da_1\, dg_0$ is a right-invariant measure on $P$, where $\Delta (a_1)=|det {\rm Ad}\, (a_1)|_{{\goth g}_0}|$. Since $A\subset P$ and $U\subset P$, we have $G=KP$ and, by [HC70], p.17, $dg=dk\,d_rp$ is a Haar measure on $G$.
\par\noindent
Define the character $\chi_0$ of $\goth g$ by 
$$\chi_0 (X)={\rm e}^{-2\pi i\langle X,\, X_0\rangle}\quad (X\in{\goth g}).$$
Then $G$ acts on $\chi_0$ by $g. \chi_0 (X)=\chi_0 ({\rm Ad}\, (g^{-1})X)\ (X\in{\goth g})$ and ${\rm Stab}\, \chi_0$ is equal to $G_0$. We will consider the representation
$$\pi ={\rm Ind}_{{\goth g}\ltimes G_0\,\uparrow\, {\goth g}\ltimes G}\,  \chi_0\otimes 1,$$
and will determine its trace. Therefore we apply formula (1). Set $G_1={\goth g}\ltimes G$. We obtain:
$${\rm tr}\,\pi\,(\varphi)=\int_{G/G_0}\int_{\goth g}\int_{G_0} \,\varphi ({\rm Ad}\ (g)X,\, gg_0g^{-1})\,\chi_0(X)\, dg_0 dX d{\dot g}$$
for $\varphi\in C_c^\infty (G_1)$. Let us determine the Haar measures more explicitly. We have
\begin{eqnarray*}
dg&=& \Delta(a_1)\, da_1\,dk\,dg_0 ,\\
d{\dot g}&=& \Delta(a_1)\,dk\, da_1.
\end{eqnarray*}
We thus obtain:
$${\rm tr}\,\pi\, (\varphi )=\int_{A_1}\int_{\goth g}\int_{G_0}\, \varphi ({\rm Ad}\, (a_1)X,\, a_1g_0a_1^{-1})\, \chi_0(X)\, \Delta(a_1)\, da_1\,dX\, dg_0$$
where we have taken $\varphi$ to be ${\rm Ad}\, (K)$-invariant. Taking $\varphi$ even of the form $\varphi (X,g)=\varphi_1(X)\,\varphi_2(g)\ (X\in{\goth g}, g\in G)$ with $\varphi_1\in C_c^\infty ({\goth g}),\, \varphi_2\in C_c^\infty (G)$ we get:
$${\rm tr}\, \pi\, (\varphi)=\int_0^\infty \widehat{\varphi}_1 (\lambda^{\Vert H_0\Vert^2}\,X_0)\,{d\lambda\over\lambda}\cdot \int_{G_0}\,\varphi_2(g_0)\, dg_0.$$
Setting $\mu =\lambda^{\Vert H_0\Vert ^2}$ we get, up to the constant ${\Vert H_0\Vert}^{-2}$
$${\rm tr}\, \pi\, (\varphi)=\int_0^\infty \widehat{\varphi}_1 (\mu\,X_0)\,{d\mu\over\mu}\cdot \int_{G_0}\,\varphi_2(g_0)\,dg_0.$$
This clearly gives a divergent first integral for suitable functions $\varphi_1$. Thus $G$ must be compact. The converse follows from Corollary 7.2 and the observation that $\goth g\ltimes G$ is a Mackey-regular semi-direct product for $G$ compact. See e.g. [Kla79], p. 105.\qed
\vskip.2cm
Though the group $\goth g\ltimes G$ is not trace class if $G$ is noncompact, it is still, at least if $G$ is algebraic, a type I group by a result of Dixmier [D57a]. Therefore sufficiently many trace class representations exist to compose the Plancherel formula. The next section is devoted to this topic. We start with complex groups. The general case requires more technical details.
\section{Plancherel formula}
Let $G$ be a complex, semisimple, connected Lie group. Then $G$ is algebraic and has finite center. In particular $G$ is a type I group. Let $\goth g$ be the Lie algebra of $G$. Select a compact real form $\goth u$ of $\goth g$ and let $U$ be the compact real analytic subgroup of $G$ with Lie algebra $\goth u$.
\par
Select, in addition, a maximal abelian subspace $\goth a$ of $i\goth u$ and set $\goth t=i\goth a +\goth a$. Then $\goth t$ is the Lie algebra of a complex maximal torus $T$, $T={\rm exp}\, \goth t ={\rm exp}\, i\goth a\,\cdot\, {\rm exp}\, \goth a$. Observe that ${\rm exp}\,\goth i\goth a$ is compact. The elements of the dual group $\widehat T$ are given by the the unitary characters of the form $\rho_{n,s}$ with $n\in L$, $L\subset {i\goth a}^\ast$ being the weight lattice of $(\goth u,\, i\goth a)$, and $s\in {\goth a}^\ast$,
$$\rho_{n,s}(t)=\rho_{n,s}({\rm exp}\, H)={\rm e}^{-2\pi i n(H_1)}\cdot {\rm e}^{-2\pi i s(H_2)},$$
where $H=H_1 +H_2, H_1\in\goth i\goth a\,,\, H_2\in\goth a$. 
\par
For $X\in \goth g $ set 
$${\rm det}\, (z - {\rm ad}\, X)=\eta(X)\, z^l + \cdots {\rm(terms\ of\ higher\ degree\ in}\ z),$$
where $z$ is an indeterminate  and $l$ is the complex rank of $\goth g$. The element $X$ is called {\it regular} if $\eta (X)\not= 0$. The set of regular elements of $\goth g$ is usually denoted by ${\goth g}'$. 
\par
Denote by  $\langle\ ,\ \rangle$ the Killing form of $\goth g$. Select an additve unitary character of the real vector space $\goth g$ of the form
$$\chi_0(X)={\rm e}^{-2\pi i{\rm Re}\, \langle X,H_0\rangle}\quad (X\in\goth g),$$
with $H_0$ a regular element in $\goth t$. Notice that $Z_G (H_0)=T$. Set
$$\pi=\pi_{H_0,n,s}={\rm Ind}_{\goth g\ltimes T\uparrow \goth g\ltimes G}\, \chi_0\otimes \rho_{n,s}.$$
We shall prove that ${\rm tr}\, \pi$ exists and then show that the representations $\pi_{H_0,n,s}$ can be used to compose the Plancherel formula. We start again with formula (1):
\begin{equation}
{\rm tr}\ \pi(\varphi)=\int_{G/T}\int_{\goth g}\int_T \,\varphi ({\rm Ad}\, (g)X,gtg^{-1})\, \chi_o(X)\,\rho_{n,s}(t) \, dt dXd\dot{g}
\end{equation}
for $\varphi\in C_c^\infty (\goth g\ltimes G).$ This can also be written as:
\begin{equation}
{\rm tr}\ \pi(\varphi)=\int_{G/T}\int_T \,\widehat\varphi ({\rm Ad}\, (g)H_0,gtg^{-1})\,\rho_{n,s}(t) \, dt d\dot{g},
\end{equation}
where $\widehat\varphi$ is the Fourier transform of $\varphi$ with respect to the first argument. 
Set $A={\rm exp}\, \goth a$ and let $G=UNA$ be an Iwasawa decomposition of $G$. Then $G/T$ can be identified with $U/{\rm exp}i\goth a\,\cdot N$ and $d\dot{g}=d\dot{u}\, dn$. Form (4) we thus obtain
\begin{equation}
{\rm tr}\ \pi(\varphi)=\int_N\int_T \,\widehat\varphi ({\rm Ad}\, (n)H_0,ntn^{-1})\, \rho_{n,s}(t) \, dt dn,
\end{equation}
where $\varphi$ is assumed to be ${\rm Ad}\, (U)$-invariant. Writing
$$ntn^{-1}=t\cdot (t^{-1}ntn^{-1})= tn'$$
with $n'\in N$, we see that $t$ varies in a compact set. Furthermore, similarly to [He67, Ch. IV, Lemma 4.6] and [Wa73, 7.5.15], we have
$$\int_N \widehat\psi ({\rm Ad}(n)\, H_0)\, dn=\int_{\goth n} \widehat\psi (H_0 +X)\, |\eta\,(H_0)|^{-2}\, dX$$ 
for $\psi\in C_c^\infty (\goth g)$, where  $\goth n$ is the Lie algebra of $N$ and $dX$ the Lebesgue measure on $\goth n$ corresponding to $dn$ by the exponential mapping.These two facts together show that ${\rm tr}\, \pi(\varphi)$ exists. Let us now proceed to determine the Plancherel formula. A main ingredient will be the following integration formula ([Wa73, 7.8.3]):
\begin{equation}
\int_{\goth g}\psi (X)\, dX={1\over |W_T|}\int_{\goth t}\int_{G/T}|\eta (H)|^2\, \psi ({\rm Ad}\, (g) H)\, d\dot{g}dH,
\end{equation}
for $\psi\in L^1(\goth g)$. In the constant $|W_T|$, $W_T$ stands for the group $N_G(T)/Z_G(T)$. From (4) we obtain:
$$\int_{G/T} \widehat\varphi ({\rm Ad}\, (g)\, H_0,1)\, d\dot{g}=\sum_n \int_{{\goth a}^\ast}\, {\rm tr}\,\pi_{H_0,n.s}(\varphi)\, ds.$$
Multiplying both sides with ${1\over |W_T|}\, |\eta(H_0)|^2$ and integrating over $\goth t$ we obtain:
\setcounter{stelling}{0}
\begin{stelling} {\rm (Plancherel formula)}
$$\varphi (0,1)={1\over |W_T|}\int_{\goth t}\int_{{\goth a}^\ast} \sum_{n\in L} {\rm tr}\,\pi_{H_0,n,s}(\varphi)\, |\eta (H_0)|^2\, ds\, dH_0$$
for any $\varphi\in C_c^\infty(\goth g\ltimes G)$.
\end{stelling}
\section{Plancherel formula: the general case}
Let $\bf G$ be a connected, complex, semisimple algebraic group defined over $\Bbb R$ and let $G={\bf G}(\Bbb R)$ denote the group of real points of $\bf G$. Call $\bf g$ the Lie algebra of $\bf G$ and $\goth g$ the Lie algebra of $G$. It is known that $\goth g$ is semisimple, that $G$ has finite center and has finitely many connected components in the ordinary topology.
\par
Let $\theta$ be a Cartan involution of $\goth g$, $\widetilde\theta$ its extension to $\bf g$:
$${\widetilde\theta} (X+iY)=\theta (X) - i\, \theta(Y)\quad (X,Y\in{\goth g}).$$
Denote by $\goth u\subset{\bf g}$ the real Lie algebra of fixed points of $\widetilde{\theta}$ and by $U$ the real analytic subgroup of $\bf G$ corresponding to $\goth u$. Then $U$ is compact. Set $K=U\cap G$. 
\par
Let $\bf T$ be a maximal torus in $\bf G$, defined over $\Bbb R$, that is a connected maximal abelian subgroup of $\bf G$, defined over $\Bbb R$. Set $T={\bf T}(\Bbb R)$ for the group of real points of $\bf T$. We shall say that $T$ is a maximal torus in $G$. Furthermore, denote by $\bf t$ the Lie algebra of $\bf T$, by $\goth t$ the Lie algebra of $T$. It is known that any $T$ is conjugate with respect to $G$ to a $\theta$-invariant one. Let us assume from now one that all $\goth t$ are $\theta$-invariant. It is also known that there are only finitely many non-conjugate $\theta$-invariant $\goth t$. Fix a maximal torus $T$ with a $\theta$-invariant Lie algebra $\goth t$ and let $A$ denote its maximal split torus. Then $A$ is the split component of a parabolic subgroup $P=MAN$ of $G$ ([HC70], Lemma 18). One also has $G=KP$.
\begin{hulpst} 
The Haar measure $dg$ on $G$ can be normalized  so that if $dk$ is the normalized invariant measure on $K$, $da$ an invariant measure on $A$, $dm$ one on $M$ and $dn$ an invariant measure on $N$ and if $f\in C_c(G)$:
$$\int_G f(g)\,dg= \int_{K\times M\times A\times N}f(kman)a^{2\rho}dkdmdadn,$$
where $a^{2\rho}=|{\rm det}Ad (a)|_{\goth n}|\ (a\in A)$, $\goth n$ being the Lie algebra of $N$.
\end{hulpst}
We refer to [Wa73], 7.6.4 for a similar result. For $X\in\goth g$ set
$${\rm det}\, (t - {\rm ad} X)=\eta (X)t^l + \ldots\ ({\rm terms\ of\ higher\ degree\ in}\ t),$$
where $t$ is an indeterminate  and $l={\rm rank}\,\goth g$. The element $X$ is called {\it regular} if $\eta (X)\not= 0$. The set of regular elements of $\goth g$ is denoted by $\goth g'$.
\par
 Let $\goth t$ be as above. Set $\goth t'$ for the set of regular elements of $\goth t$ and define ${\goth g'}_T= {\rm Ad}\,(G)\, \goth t'$. Then
$$\goth g'=\bigcup_T {\goth g'}_T.$$
Observe that this is actually a finite union and we can take representatives $T$ with Lie algebra $\goth t$ being $\theta$-invariant.
\par
Denote by $\langle\ ,\ \rangle$ the Killing form of $\goth g$. Observe that $-\langle X,\, \theta(Y)\rangle\ (X,Y\in\goth g)$ is a scalar product on $\goth g$. Choose Euclidean measures $dX$ on $\goth g$ and $dH$ on $\goth t$ and normalize the Haar measures $dg$ on $G$ and $dt$ on $T$ so that at the identity they correspond to $dX$ and $dH$ via the exponential mapping. Let $d{\dot g}$ denote the invariant measure on $G/T$ such that $dg=d{\dot g} dt$. Then one has the following integration formula:
\setcounter{stelling}{2}
\begin{hulpst}
$$\int_{{\goth g'}_T}\, f(X) dX= {1\over |{W_T}|}\int_{\goth t}\, |\eta (H)|\,\int_{G/T}f({\rm Ad}(g)H)\, d{\dot g}dH\quad (f\in L^1(\goth g)).$$
where $|W_T|$ is the number of elements of $W_T=\widetilde T/T$, $\widetilde T$ being the normalizer of $T$ in $G$. 
\end{hulpst}
We refer to [Wa73], 7.8.3 and its proof. Observe that we might be dealing with another normalization of $dg$ in Lemma 11.1. 
\par
For any $T$ let $\rho^T$ denote a continuous, unitary, character of $T$ and let $d\rho^T$ denote the Haar measure on the dual group $\widehat T$, dual to the Haar measure $dt$ on $T$.
\par
Select, in addition, an additive character of the real vector space $\goth g$ of the form
$$\chi_0^T(X)={\rm e}^{-2\pi i  \langle X,\, H_0\rangle}\quad (X\in{\goth g}),$$
with $H_0$ a regular element of $\goth t$, the Lie algebra of $T$. Notice that the centralizer of $H_0$ in $G$ equals $T$. Set 
$$\pi^T=\pi^T_{H_0,\rho^T}={\rm Ind}_{\goth g\ltimes T\uparrow{\goth g}\ltimes G}\,\chi_0^T\otimes \rho^T.$$
We shall first show that ${\rm tr}\, \pi^T$ exists and then show that these representations are sufficient to compose the Plancherel formula.
\par
According to formula (1) we have:
\begin{equation}
{\rm tr}\, \pi^T(\varphi)=\int_{G/T}\int_{\goth g}\int_T\, \varphi ({\rm Ad}\, (g)X, gtg^{-1})\,\chi_0^T(X)\, \rho^T(t)\, dtdXd{\dot g},
\end{equation}
for all $\varphi\in C_c^\infty (\goth g\ltimes G)$. This can also be written as:
\begin{equation}
{\rm tr}\, \pi^T(\varphi)=\int_{G/T}\int_T\, \widehat\varphi ({\rm Ad}\, (g)H_0, gtg^{-1})\,\rho^T(t)\, dtd{\dot g},
\end{equation}
where $\widehat\varphi$ is the Fourier transform of $\varphi$ with respect to the first argument.
\par
To prove that ${\rm tr}\,\pi^T(\varphi)$ exists we apply the decomposition $G=KP=KNMA$ and $dg=dkdndmda$. As observed before, we might be using here a different normalization of $dg$, but it is clear that the exsistence of ${\rm tr}\,\pi^T(\varphi)$ does not depend on the normalization used.
\par
Let us write $T=AB$ with $B=T\cap K$ compact. Then $B\subset M\cap K$ and $G/T$ can be identified with $KN{\dot M}$ where ${\dot M}\simeq M/B$. Furthermore $d{\dot g}$ can be identified with $dkdnd{\dot m}$. Invoking this in formula (8) we obtain:
\begin{equation}
{\rm tr}\,\pi^T(\varphi)=\int_N\int_{\dot M}\int_T\widehat\varphi ({\rm Ad}\, (nm)H_0, nmtm^{-1}n^{-1})\, \rho^T(t)\, dtd{\dot m}dn,
\end{equation}
where $\varphi$ is assumed to be ${\rm Ad}(K)$-invariant. 
\par
Writing $t=ab\ (a\in A, b\in B)$ and 
$$nmtm^{-1}n^{-1}=nmabm^{-1}n^{-1}=a\cdot a^{-1}na m' n^{-1}m'^{-1} \, m'=a\cdot n'\cdot m'$$
with $m'=mbm^{-1}$, we see that $a$, and thus $t$, varies in a compact set. We also know, by  a result of Harish-Chandra, that the the invariant measure $d{\dot g}$ on ${\rm Ad}(G)H_0$ extends to a tempered measure on $\goth g$ (cf. [Var77], page 40, Corollary 10). These two facts together imply that ${\rm tr}\,\pi^T(\varphi )$ exists. 
\par
From equation (8) we now obtain:
\begin{equation}
\int_{G/T}\widehat\varphi ({\rm Ad}\,(g)H_0, 1)\,d{\dot g}= \int_{\widehat T}\, {\rm tr}_{H_0,\rho^T}(\varphi)\, d\rho^T.
\end{equation}
Multiplying both sides with ${1\over |W_T|}\, |\eta(H_0)|$ and integrating over $\goth t$ we obtain: 
\begin{equation}
\int_{\goth g'_T}\widehat\varphi(X)\, dX={1\over |W_T|}\int_{\widehat T}\int_{\goth t}\,{\rm tr}\,\pi_{H_0,\rho^T}(\varphi)\, |\eta(H_0)|\, d\rho^TdH_0,
\end{equation}
for any $\varphi\in C_c^\infty (\goth g\ltimes G)$. Finally we obtain
\begin{stelling} {\rm (Plancherel\ formula)}\par
For any $\varphi\in C_c^\infty (\goth g\ltimes G)$ one has
$$\varphi (0,1)=\sum_T\,{1\over |W_T|}\int_{\widehat T}\int_{\goth t}\,{\rm tr}\,\pi_{H_0,\rho^T}(\varphi)\, |\eta(H_0)|\, d\rho^TdH_0,$$
where the summation is over a complete (finite) set of non-conjugate maximal tori $T$ with $\theta$-invariant Lie algebra.
\end{stelling}
\section{The special case $G={\rm SL}(2,\Bbb R)$}
\setcounter{stelling}{0}
In this section we write down the explicit form of the Plancherel formula for $\goth g\ltimes G$ with $G={\rm SL}(2,\Bbb R)$. 
\par\noindent
Consider the following two tori in the group $G={\rm SL}(2,\Bbb R)$:
$$A=\left \{\left (\begin{array}{cc}a&0\\0&a\end{array}\right ):\, a\not= 0\right \}\quad {\rm and}\quad B
=\left \{\left (\begin{array}{cc}{\rm cos}\varphi&{\rm sin}\varphi\\-{\rm sin}\varphi &{\rm cos}\varphi\end{array}\right ):\, 0\leq\varphi<2\pi\right \}.$$
Let us define characters of $A$ by
$$\rho^A_{\varepsilon, s}\, \left (\begin{array}{cc}a&0\\0&a\end{array}\right )=|a|^{-2\pi i s}\,\left ({a\over |a|}\right )^\varepsilon$$
for $s\in \Bbb R,\, \varepsilon =0,1$, and characters of $B$ by 
$$\rho^B_n\left (\begin{array}{cc}{\rm cos}\varphi&{\rm sin}\varphi\\-{\rm sin}\varphi &{\rm cos}\varphi\end{array}\right )={\rm e}^{-2\pi i n\varphi}$$
for $n\in\Bbb Z$.
\par
Define the following elements of the Lie algebra $\goth a$ of $A$ and the Lie algebra $\goth b$ of $B$ respectively:
$$H_u=\left (\begin{array}{cc}u&0\\0&-u\end{array}\right )\in\goth a,\quad H_v =\left (\begin{array}{cc}0&v\\-v&0\end{array}\right )\in\goth b,$$
and take $u\not= 0,v\not= 0$.
\par
They give rise to the following additive characters of $\goth g$:
$$\chi_u(X)={\rm e}^{-2\pi i {\rm tr}\, (XH_u)},\quad \chi_v(X)={\rm e}^{-2\pi i {\rm tr}\, (XH_v)}\ \ (X\in\goth g).$$
Let us now consider the unitary representations:
$$\pi^A_{\varepsilon,s,u}={\rm Ind}_{\goth g\ltimes A\uparrow\goth g\ltimes G}\ \chi_u\otimes\rho^A_{\varepsilon, s}\quad{\rm and}\quad \pi^B_{n,v}={\rm Ind}_{\goth g\ltimes B\uparrow\goth g\ltimes G}\,\chi_v\otimes\rho^B_n.$$
It can be shown that both sets of representations are irreducible and of trace class, and that the following theorem holds:
\begin{stelling} {\rm (Plancherel theorem)}
\par\noindent
For any $\varphi\in C_c^\infty (\goth g\ltimes G)$ one has
$$\varphi (0,e)=2\sum_{\varepsilon =0,1}\int_{-\infty}^\infty\int_{-\infty}^\infty {\rm tr}\,\pi^A_{\varepsilon, s,u}(\varphi )\, |u|\, dsdu +4\sum_{n\in\Bbb Z}\int_{-\infty}^\infty {\rm tr}\,\pi^B_{n,v}(\varphi )\, |v|\, dv.$$
\end{stelling}
\section{A generalization}
\setcounter{stelling}{0}
In this section formulate a theorem which generalizes Theorem 9.1 in case of algebraic groups. Let us first fix some notation. 
\par
Let $\mathbf G$ be a connected, complex, semisimple, linear algebraic group defined over $\Bbb R$ and let $G={\mathbf G}(\Bbb R)$ its group of real points. Then $G$ is a semisimple Lie group with finite center and finitely many connected components. Let $V$ be a real finite-dimensional vector space and $\mathbf V$   its compexification. Assume that $\mathbf G$ acts on $\mathbf V$ by a representation $\rho$ which is, in addition, an algebraic morphism defined over $\Bbb R$. Set $G_1=V\ltimes G$. Let $\mathbf H$ be the kernel of $\rho$ and set $H={\mathbf H}(\Bbb R)$.
\begin{stelling} $G_1$ is a trace class group if and only if $G/H$ is compact.
\end{stelling} 
As to the proof of this theorem, we can only confirm the `if'- part at the moment. We formulate it in the following lemma.
\setcounter{stelling}{2}
\begin{hulpst} If $G/H$ is compact then the sem-idirect product $G_1=V\ltimes G$ is a trace class group.
\end{hulpst}
\par\noindent
{\bf Proof}. First of all observe that that the orbit space $V^\ast\backslash G$ is equal to the space $V^\ast\backslash (G/H)$. Since $G/H$ is compact we may conclude that $V\ltimes G$ is a Mackey-regular semi-direct product. See e.g. [Kla79], p. 105.
\par
Observe that $H$ is a normal subgroup of $G$, hence semisimple with finitely many connected components, hence unimodular and a trace class group by [DD15], Proposition 1.10 (ii).  Now apply [DD15], Theorem 4.1 with $Q=G$ and $Q_0=H$. With the notation of this theorem, observe that for any unitary character $\chi$ of $V$ its stability group $Q_\chi$ in $Q$
is unimodular and trace class, because of [KL72], p. 470. Hence $G_1$ is a trace class group. 
 \qed

\vskip.2cm
Gerrit van Dijk
\par\noindent
Mathematisch Instituut
\par\noindent
Niels Bohrweg 1
\par\noindent
2333 CA Leiden, The Netherlands
\par\noindent
dijk@math.leidenuniv.nl
\end{document}